\documentclass[12pt]{amsart}
\usepackage{amssymb}

\newcommand{\Z}{{\mathbf Z}}
\newcommand{\G}{{\mathbf G}}
\newcommand{\C}{{\mathbf C}}

\newcommand{\Q}{{\mathbf Q}}
\newcommand{\bP}{{\mathbf P}}
\newcommand{\cO}{{\mathcal O}}
\DeclareMathOperator{\cl}{cl}
\DeclareMathOperator{\image}{Im}
\DeclareMathOperator{\chara}{char}
\DeclareMathOperator{\Spec}{Spec}
\DeclareMathOperator{\Pic}{Pic}
\DeclareMathOperator{\Br}{Br}
\theoremstyle{plain}
\newtheorem  {thm}        {Theorem}
\newtheorem  {lemma}      {Lemma}
\newtheorem  {prop}       {Proposition}
\newtheorem  {cor}        {Corollary}
\theoremstyle{definition}

\theoremstyle{remark}

\begin{document}

\title[Hodge-theoretic obstruction]{
  Hodge-theoretic obstruction to existence of quaternion algebras
}
\author{Andrew Kresch}
\address{
  Department of Mathematics,
  University of Pennsylvania,
  Philadelphia, PA 19104
}
\email{kresch@math.upenn.edu}
\date{September 12, 2000}
\thanks{The author was partially supported by
an NSF Postdoctoral Research Fellowship.}

\maketitle

\section{Introduction} \label{s.intro}

The subject of this paper is the Brauer group of a nonsingular complex
projective variety.
More specifically, we study the question of whether a 2-torsion element
of the cohomological Brauer group is representable by a quaternion algebra
over the generic point.
Using intersection theory -- on schemes and on algebraic stacks -- we are able
to describe an obstruction to such a representation, and therefore, to give
examples of varieties with 2-torsion classes that are
not representable by quaternion algebras.

Let $X$ be a nonsingular complex algebraic variety.
It is a well-known consequence of the Grauert--Remmert theorem, plus GAGA,
that topological 2-sheeted covers,
classified by $H^1(X,\Z/2)$, correspond exactly to algebraic
degree 2 unramified covers.
There is similarly an algebraic object that we can associate to an
element $\alpha\in H^2(X,\Z/2)$ -- a kind of algebraic stack over $X$ known as
a {\em gerbe}.
The element $\alpha$ determines a
2-torsion element of $\Br(X)=H^2(X,\G_m)$, and if
this element is represented by a quaternion algebra then the quaternion algebra
can be identified with the (descent to $X$ of the) endomorphism algebra of a
rank 2 vector bundle on the gerbe.
The second Chern class of this bundle (or rather, a specific multiple of this)
is an algebraic class on $X$, and if $X$ is projective this leads to a
Hodge-theoretic obstruction to representing the element
of $\Br(X)$ by a quaternion algebra.

We turn now to the statement of our main result.

\begin{thm}
\label{thm.main}
Let $X$ be a nonsingular complex projective variety.
Let $\alpha\in H^2(X,\Z/2)$ be a class such that
{\em (i)} the image of $\alpha$ in the Brauer group of the function field
$\Br(k(X))$ is representable by a quaternion algebra over $k(X)$, and
{\em (ii)} there exists a preimage $\alpha_0\in H^2(X)$ under the
natural map $H^2(X)=H^2(X,\Z)\to H^2(X,\Z/2)$.
Then there exists an algebraic class $\eta\in A^2X$ whose cycle class
in cohomology satisfies
$$\cl(\eta)=4(\alpha_0^2+2\varepsilon)$$
for some $\varepsilon\in H^4(X)$.
\end{thm}

An immediate consequence of Theorem \ref{thm.main} is the following result.

\begin{cor}
\label{cor.hodge}
Let $X$ be a nonsingular complex projective variety.
Let $N$ denote the Neron--Severi group of $X$, and
$i\colon (H^2(X)/N)\otimes\Z/2\to \Br(k(X))$ the injective
homomorphism to the Brauer group induced by reduction of
coefficients $H^2(X)\to H^2(X,\Z/2)$, followed by the natural map to
the Brauer group.
Let $M=H^4(X)\cap H^{2,2}(X)$; then for any
$\beta\in (H^2(X)/N)\otimes\Z/2$,
the image of $\beta$ under the homomorphism
\begin{equation}
\label{obstr}
(H^2(X)/N)\otimes\Z/2\to (H^4(X)/M)\otimes\Z/2
\end{equation}
given by $\beta\mapsto\beta^2$
is an obstruction to representing
$i(\beta)$ by a quaternion algebra over $k(X)$.
\end{cor}

After stating conventions and preliminary results in Section \ref{s.prelim},
we prove Theorem \ref{thm.main} in Section \ref{s.proof}.
In Section \ref{s.example}, we provide examples for which the
obstruction map (\ref{obstr}) is nontrivial.
The geometric significance of a nontrivially obstructed class $\beta$ is that
$i(\beta)$ cannot be represented by a Brauer--Severi variety of
relative dimension 1, even over the generic point of $X$.

A deep result of Merkurev \cite{Me} asserts that
if $F$ is any field with $\chara F\ne 2$,
then the 2-torsion subgroup of the Brauer group is generated by the classes
of quaternion algebras.
Then we say $F$ is {\em linked} if every 2-torsion element of
$\Br(F)$ is representable by a quaternion algebra.
Some classes of fields known to be linked are listed in \cite{E-L}.
This list includes all fields of transcendence degree 2 over $\C$;
it is shown in \cite{Sar} that more generally any $C_2$-field is linked.
While \cite{E-L} gives some explicit examples of non-linked fields,
the important main result in that paper there is
a necessary (but not sufficient) characterization of linked fields by
their $u$-invariant, and this provides a way to generate unlimited
numbers of examples of non-linked fields.

In Section \ref{s.example} we present a result, Theorem \ref{thm.example},
which shows that
function fields of complex algebraic threefolds are not, in general, linked.
The threefolds we study are
Brauer--Severi varieties over surfaces, and for such $X$ there are
classes in $H^2(X,\Z/2)$ for which the induced elements of
$\Br(k(X))$ are not representable by quaternion algebras over $k(X)$.
So, our examples differ from previously known examples in that the
constructions are of a geometric nature.

As a consequence of cohomological purity and the valuative criterion
for properness, the image of the inclusion
$\Br(X)\to\Br(k(X))$ can be identified with the
unramified Brauer group \cite{Sal}
of $k(X)$.
Thus, the cohomological Brauer group is a birational invariant,
for nonsingular complete
varieties over a given field of characteristic zero.
It is therefore a weakness of the present approach to require a specific
nonsingular projective model to study a problem which is intrinsic
to the function field.

The author would like to thank Brendan Hassett and Raymond Hoobler for helpful
correspondence.

\section{Preliminaries} \label{s.prelim}

All schemes are of finite type over a field.
All stacks are algebraic and of finite type over a field.
$A_*$ denotes the Chow groups of a stack as in \cite{Kr}.
For a smooth stack $X$ of pure dimension $n$, we let $A^kX$ denote
$A_{n-k}X$.
For a topological space $X$, $H^*(X)$ denotes $H^*(X,\Z)$.
When $X$ is a complex variety (respectively, a reduced algebraic stack of
finite type over $\Spec\C$),
$H^*(X)$ denotes the cohomology of the underlying
analytic space (respectively, a topological realization of the underlying
analytic stack).
Relevant facts on stacks can be found in \cite{D-M,L-MB}.
$\Br(X)$ denotes the cohomological Brauer group of $X$;
facts on Brauer groups can be found in \cite{Gr,Milne}.

\begin{prop}
\label{p.pica1}
Let $X$ be a regular algebraic stack of pure dimension $n$,
and assume $X$ contains an open
substack $U$ such that {\em (i)} the complement $Z=X\smallsetminus U$ is
empty or of codimension $\ge 2$ in $X$, and {\em (ii)} $U$ is isomorphic
to a quotient stack of the form $[U/G]$, where $U$ is an algebraic space
and $G$ is a linear algebraic group.
Then the first Chern class homomorphism $\Pic(X)\to A_{n-1}X$ is an isomorphism.
\end{prop}

\begin{proof}
There are obvious isomorphisms $\Pic(X)\to \Pic(U)$ and
$A_{n-1}X\to A_{n-1}U$, compatible with the first Chern class map, so
it suffices to consider the case $X=U$ is a quotient stack,
and we are reduced to \cite[Thm.\ 1]{E-G}.
\end{proof}

\begin{cor}
\label{cor.pica1}
Let $X$ be a regular algebraic stack of pure dimension $n$ with finite
stabilizers at all geometric points.
Then $\Pic(X)\simeq A_{n-1}X$.
\end{cor}

\begin{proof}
A finite cover by a scheme $Y\to X$ exists by \cite[Thm.\ 2.8]{FourAuthors}
(or in the case of primary interest -- $X$ a Deligne-Mumford stack --
by \cite[Thm.\ 16.6]{L-MB}).
We may suppose $Y$ normal, hence finite flat over $X$ away from
a closed substack $Z$ of $X$ which is empty or of codimension $\ge 3$.
This implies (cf.\ \cite[Prop.\ 3.5.7]{Kr}) that
$U:=X\smallsetminus Z$ is a quotient stack,
and Proposition \ref{p.pica1} applies.
\end{proof}

Results that follow use the fact that gerbes over a scheme $X$, banded
by the group of $n^{\rm th}$ roots of unity $\mu_n$,
are classified by $H^2(X,\mu_n)$
\cite[\S IV.2]{Milne}.
When $n$ is invertible in the base field (the case of primary
interest to us), this classifying group is an \'etale cohomology group;
in general, cohomology for the flat (fppf) topology must be employed.

\begin{prop}
\label{p.tfae}
Let $X$ be a regular scheme, and let $n$ be a positive integer.
Let $\beta\in H^2(X,\mu_n)$, and let ${\mathcal G}$ be the
gerbe over $X$, banded by $\mu_n$, classified by $\beta$.
Then the following are equivalent:
\begin{itemize}
\item[(i)] The gerbe ${\mathcal G}\to X$ is Zariski locally trivial.
\item[(ii)] ${\mathcal G}$ is a trivial gerbe over some nonempty Zariski
open subset of $X$.
\item[(iii)] The image of $\beta$ in $\Br(k(X))$ is zero.
\item[(iv)] There exists a line bundle on ${\mathcal G}$ on which
the action of stabilizer groups at geometric points of ${\mathcal G}$
is faithful.
\end{itemize}
\end{prop}

\begin{proof}
Clearly, (i) implies (ii), and (ii) implies (iii).
As $\Br(X)\to \Br(k(X))$ is injective, (iii) implies
$\beta$ lies in the image of the boundary homomorphism
$$\delta\colon H^1(X,\G_m)\to H^2(X,\mu_n)$$
of the Kummer sequence.
If $\beta=\delta(\alpha)$, then ${\mathcal G}$ is isomorphic to a
$\G_m$-quotient of the principal bundle on $X$ associated to $\alpha$.
So, (iii) implies (iv).
Given a line bundle as in (iv), we can identify ${\mathcal G}$ with
a $\G_m$-quotient of the associated principal bundle $P$.
Now $P\to X$ is the principal bundle of a class $\alpha\in H^1(X,\G_m)$
with $\delta(k\alpha)=\beta$ for some $k$ prime to $n$, and (i) holds.
\end{proof}

\begin{cor}
\label{cor.picgerbe}
Let $p$ be a prime.
Let $X$ be a regular scheme, $\beta\in H^2(X,\mu_p)$, and let
$f\colon {\mathcal G}\to X$ be the gerbe banded by $\mu_p$ with class $\beta$.
Then
$$f^*\colon\Pic(X)\to \Pic({\mathcal G})$$
is an isomorphism if and only if the image of $\beta$ in $\Br(k(X))$
is nonzero.
\end{cor}

\begin{prop}
\label{p.algebraicgerbe}
Let $n$ be a positive integer.
Let $X$ be a scheme, and let ${\mathcal G}$ be the gerbe over $X$,
banded by $\mu_n$, classified by $\beta\in H^2(X,\mu_n)$.
Now suppose $\beta$ is in
the image under the boundary homomorphism of the Kummer sequence.
If $E$ is a line bundle on $X$ whose class in $\Pic(X)$ maps under
the boundary homomorphism to $\beta$, then
for each $k$ there is an exact sequence
$$\bigoplus_{s\ge 1} A_{k+s}X\to
\bigoplus_{s\ge 0}A_{k+s}X\to A_k{\mathcal G}\to 0$$
where the map on the left is
$$(\alpha_1,\alpha_2,\ldots)\mapsto
(c_1(E)\cap\alpha_1, c_1(E)\cap\alpha_2-n\alpha_1,\ldots).$$
\end{prop}

\begin{proof}
We follow the program of Equivariant Intersection Theory \cite{E-G}:
${\mathcal G}$ is a $\G_m$-quotient of $E$ minus the zero section,
so scheme approximations
to ${\mathcal G}$ are complements of zero sections of line bundles over
$X\times\bP^r$:
$$A_k{\mathcal G}=
A_{k+r+1}((E\otimes\cO_{\bP^r}(-n)\smallsetminus s(X\times\bP^r))$$
for $r$ sufficiently large.
The result now follows from the excision sequence
$$A_{k+r+1}(X\times\bP^r)\to A_{k+r+1}(E\otimes\cO_{\bP^r}(-n))\to
A_k{\mathcal G}\to 0.\qed$$
\renewcommand{\qed}{}\end{proof}

\begin{cor}
\label{cor.smoothalgebraicgerbe}
If the conditions of Proposition \ref{p.algebraicgerbe} are met and
additionally $X$ is smooth, then
there is a natural map
$$A^*(X)[z]/(nz-c_1(E))\to A^*{\mathcal G},$$
which is an isomorphism of rings.
\end{cor}

An analogous result to Corollary \ref{cor.smoothalgebraicgerbe}
in the topological category is Proposition \ref{p.topologicalgerbe},
which we give below after a topological lemma.

\begin{lemma}
\label{l.topological}
Let $X$ be a paracompact topological space,
and let $E$ be a topological complex line bundle on $X$,
with zero section $s$.
Then there is a natural injective homomorphism for every $k$
$$H^k(X)/(c_1(E) \cup H^{k-2}(X))\to H^k(E\smallsetminus s(X)),$$
and this is an isomorphism precisely when
$c_1(E)\cup{-}\colon H^{k-1}(X)\to H^{k+1}(X)$
is injective.
\end{lemma}

\begin{proof}
We compactify $E$ by setting
$P={\mathbf P}(E\oplus 1)$, the (complex) projectivization of
the Whitney sum of $E$ and a trivial complex line bundle.
Now we have $E$ and $F:=P\smallsetminus s(X)$, complex line bundles over $X$,
with $E\cap F = E\smallsetminus s(X)$.
The Mayer-Vietoris sequence in cohomology gives
\begin{equation}
\label{mvseq}
H^k(P)\to \bigl(H^k(X)\bigr)^2\to H^k(E\smallsetminus s(X))\to
H^{k+1}(P)\to \bigl(H^{k+1}(X)\bigr)^2.
\end{equation}
Let $L$ denote the tautological complex line bundle on $P$;
then the projective bundle theorem dictates
$H^k(X)\oplus H^{k-2}(X)\simeq H^k(P)$ via
$(\alpha,\beta)\mapsto \alpha+(c_1(L)\cup\beta)$.
Since $L$ is trivial on $E$ and $L|_{P\smallsetminus E}\simeq E$,
the leftmost map in (\ref{mvseq}) is
$$(\alpha,\beta)\mapsto\bigl(\alpha,\alpha+(c_1(E)\cup\beta)\bigr).$$
The result is now immediate.
\end{proof}

Fix now a realization of the classifying space $B(\Z/n)$
as a topological abelian group \cite{Milnor}.
Then principal $B(\Z/n)$-bundles are classified by
$H^2({-},\Z/n)$.

\begin{prop}
\label{p.topologicalgerbe}
Let $X$ be a paracompact topological space, and assume $X$ is
homotopy equivalent to a CW complex.
Suppose $\beta\in H^2(X,\Z/n)$ is the image of $\beta_0\in H^2(X)$.
Let $G\to X$ be a principal $B(\Z/n)$-bundle classified by $\beta$.
Then there is an injective ring homomorphism
$$H^*(X)[u]/(nu-\beta_0)\to H^*(G),$$
which is an isomorphism in degrees $\le 2$.
\end{prop}

\begin{proof}
Let $\cO_{\C\bP^\infty}(-n)$ denote the $n^{\rm th}$ twist of
the tautological line bundle on $\C\bP^\infty$.
The complement of the zero section in $\cO_{\C\bP^\infty}(-n)$ is
homotopy equivalent
to $B(\Z/n)$, so by standard arguments, if $L$ is a complex
line bundle on $X$ with $c_1(L)=\beta_0$, then the complement, over
$X\times\C\bP^\infty$, of the zero section of $L\otimes \cO_{\C\bP^\infty}(-n)$
is homotopy equivalent to $G$.
Now the result follows from Lemma \ref{l.topological},
plus the general fact that $H^0(X)$ and $H^1(X)$ are torsion free.
\end{proof}

\begin{prop}
\label{p.timestwo}
Let $X$ be a scheme of pure dimension $k$, let $n$ be a positive
integer, and let
$f\colon {\mathcal G}\to X$ be a gerbe banded by $\mu_n$.
Then for any $\alpha\in A_{k-1}{\mathcal G}$, we have
$n\alpha=f^*\beta$ for some $\beta\in A_{k-1}X$.
\end{prop}

\begin{proof}
We may assume $X$ is reduced.
Consider the components of the regular locus
$X^{\rm reg}=\coprod X_i$.
For each $i$, if we let $f_i$ denote the restriction of $f$
over $X_i$, then we claim
the image of $f_i^*\colon A_{k-1}X_i\to A_{k-1}{\mathcal G}_i$ contains
$nA_{k-1}{\mathcal G}_i$.
Indeed, this follows by Proposition \ref{p.algebraicgerbe} when
$f_i$ is Zariski locally trivial; otherwise by Corollary \ref{cor.picgerbe}
combined with Corollary \ref{cor.pica1}, $f_i$ induces an isomorphism
of Chow groups in dimension $k-1$.
Now the desired result follows by comparing the excision sequences for
$X^{\rm reg}\subset X$ and ${\mathcal G}^{\rm reg}\subset {\mathcal G}$.
\end{proof}

\section{Proof of the main theorem} \label{s.proof}

Let us start by recalling conditions (i) and (ii) on
$\alpha\in H^2(X,\Z/2)$
(\'etale cohomology, which by cohomology comparison,
equals singular cohomology).
Condition (i), the existence of a quaternion algebra, is equivalent to
existence of a Brauer--Severi variety of dimension $1$ representing
the class of $\alpha$ in $\Br(k(X))$.
This spreads out over a nonempty open subset of $X$ and corresponds to
a rank 2 vector bundle on an open subset of the $\Z/2$-gerbe on $X$ classified
by $\alpha$.
Since $X$ is regular, the vector bundle -- and hence also the Brauer--Severi
variety -- can be extended over $X\smallsetminus Z$
for some closed $Z\subset X$,
empty or of codimension $\ge 3$.

A Brauer--Severi variety determines a class in the Brauer group, which
is the obstruction to identification with the projectivization of an
algebraic vector bundle.
The analogous topological obstruction lies in the group $H^3(X)$.
Requiring
as in (ii) that $\alpha$ lifts to some $\alpha_0\in H^2(X)$ is equivalent
to requiring that the topological obstruction vanishes.
So the setting of the theorem is the case where the obstruction vanishes
topologically but not algebraically.

We fix some notation.
We denote by
${\mathcal G}$ the gerbe on $X$, banded by $\Z/2$, classified
by $\alpha$.
We let $V$ denote a Brauer--Severi variety over $X\smallsetminus Z$,
and $E$ a rank $2$ vector bundle on
${\mathcal G}_{X\smallsetminus Z}:={\mathcal G}\times_X(X\smallsetminus Z)$
such that $\bP(E)\simeq V\times_X{\mathcal G}$.
The restrictions of ${\mathcal G}$, $V$, and $E$ over
the generic point
$\Spec k(X)$ are denoted
${\mathcal G}_{\rm gen}$, $V_{\rm gen}$, and $E_{\rm gen}$, respectively.
We note that ${\mathcal G}_{X\smallsetminus Z}$ is a
quotient stack of a smooth variety by a linear algebraic group
(e.g., of the principal bundle associated to $E$ by the group $GL_2$),
and hence $A^*{\mathcal G}_{X\smallsetminus Z}$ and
$H^*({\mathcal G}_{X\smallsetminus Z})$ can be identified with
equivariant Chow groups and equivariant cohomology groups, respectively.

Now we proceed with the proof.
The conic $V_{\rm gen}$ is not rational.  So, the
Chow ring of $V_{\rm gen}$ is $\Z[y]/(y^2)$, with $y$
the class of a degree $2$ point on $V_{\rm gen}$.
Comparing the projective bundle formula and
Corollary \ref{cor.smoothalgebraicgerbe}, we have
$$A^*{\mathcal G}_{\rm gen}[w]/(w^2-c_1(E_{\rm gen})w+c_2(E_{\rm gen}))\simeq
\Z[y,z]/(y^2,2z-ky),$$
for some integer $k$.
We have $A^1{\mathcal G}_{\rm gen}=0$ by Corollary \ref{cor.picgerbe}.
Hence $A^1{\mathbf P}(E_{\rm gen})\simeq\Z$, so $k$ must be odd, and
without loss of generality we can suppose $k=1$.
Then $w$ and $z$ both generate $A^1{\mathbf P}(E_{\rm gen})$,
so $w=\pm z$,
and now
\begin{equation}
\label{zmod4}
A^2{\mathcal G}_{\rm gen}\simeq\Z/4,
\end{equation}
with $c_2(E_{\rm gen})$ as a generator.

Topologically, we have by Proposition \ref{p.topologicalgerbe}
an injective ring homomorphism
$$H^*(X)[u]/(2u-\alpha_0)\to H^*({\mathcal G}),$$
which is an isomorphism in degrees $\le 2$.
Moreover, the obstruction to writing $V$ as the projectivization of a
complex rank $2$ vector bundle vanishes topologically, so we may write
$$V\simeq {\mathbf P}(B)$$
where $B$ is a topological rank $2$ complex vector bundle
on $X\smallsetminus Z$.
Let $f$ denote the restriction of the
projection ${\mathcal G}\to X$, over $X\smallsetminus Z$.
We have ${\mathbf P}(f^*B)\simeq {\mathbf P}(E)$, and hence
$$f^*B\simeq E\otimes L$$
for some topological complex line bundle $L$ on
${\mathcal G}_{X\smallsetminus Z}$.
Comparing the restrictions over a point of $X$, we see that $L$ must
be nontrivial on fibers of $f$.
So,
\begin{equation}
\label{c1class}
c_1(L)=u+f^*\delta
\end{equation}
for some $\delta\in H^2(X\smallsetminus Z)\simeq H^2(X)$.
The cycle class map to equivariant co\-ho\-mology \cite{E-G} respects
Chern classes; using $A^1{\mathcal G}\simeq A^1X$
(Corollary \ref{cor.picgerbe}), we see
$$c_2(E)=-c_1(L)^2-c_1(E)c_1(L)+f^*c_2(B)=-u^2+uf^*\beta+f^*\varepsilon$$
for some $\varepsilon\in H^4(X\smallsetminus Z)\simeq H^4(X)$ and
$\beta\in H^2(X)$.  Hence
$$2c_2(E)+2u^2=f^*\beta'$$
in $H^4({\mathcal G}_{X\smallsetminus Z})$, for some $\beta'\in H^4(X)$.

Let
$\gamma=c_2(E)\in A^2{\mathcal G}_{X\smallsetminus Z}\simeq A^2{\mathcal G}.$
By (\ref{zmod4}),
$4c_2(E_{\rm gen})=0$ in $A^2{\mathcal G}_{\rm gen}$, and hence
$4\gamma$ vanishes in $A^2{\mathcal G}_U$ for some nonempty open $U\subset X$.
Let $Y=X\smallsetminus U$; we may assume $Y$ has pure dimension $n-1$, where
$n=\dim X$.
Let $i$ denote the inclusion ${\mathcal G}_Y\to {\mathcal G}$.
By excision, $4\gamma$ lies in the image of
$i_*\colon A_{n-2}{\mathcal G}_Y\to A_{n-2}{\mathcal G}$.
By Proposition \ref{p.timestwo}, now, we have
$$8\gamma = f^*\eta$$
for some $\eta\in A^2X$.
It follows that
$$f^*\cl(\eta)=-2f^*\alpha_0^2+4f^*\beta'.$$
The kernel of
$f^*\colon H^4(X)\simeq H^4(X\smallsetminus Z)
\to H^4{\mathcal G}_{X\smallsetminus Z}$
is $2$-torsion, so
$$4\alpha_0^2-8\beta'\in\image(\cl\colon A^2X\to H^4(X)),$$
and the theorem is proved.

\section{example} \label{s.example}

We show that if $X$ is any nonsingular complex projective surface
with $p_g(X)\ne 0$, then there exist Brauer--Severi varieties
$V\to X$ such that the obstruction map (\ref{obstr})
of Corollary \ref{cor.hodge} for $V$ is nontrivial.

\begin{thm}
\label{thm.example}
Let $X$ be a nonsingular complex projective surface with
nonzero geometric genus.
Let $\beta$ be an element of $H^2(X)$ whose reduction to
$H^2(X,\Z/2)$ maps to some nonzero $\lambda\in \Br(X)$.
If $V\to X$ is a smooth conic representing $\lambda$ (this exists
since $X$ is a nonsingular surface over $\C$), then the obstruction
map of Corollary \ref{cor.hodge} for the variety $V$
is nontrivial.
In fact, any element of $H^2(V)$ which generates
$H^2(V)/H^2(X)\simeq\Z$ has nonzero image under the obstruction map.
\end{thm}

\begin{proof}
Denote by $r$ the Picard number of $X$; then $r$ is strictly less
than the second Betti number $k$ of $X$.
Let $N$ be the Neron-Severi group of $X$.
Without loss of generality, we may suppose the image of $\beta$ in
$H^2(X)/N$ is a primitive lattice element.
Then we can write
$$H^2(X) = N \oplus \langle\beta\rangle \oplus T,$$
where $T$ has rank $k-r-1$.

Recall that the hypotheses dictate that $V$ is topologically
(but not algebraically)
the projectivization of a rank 2 complex vector
bundle on $X$.
To compute the obstruction map for $V$,
we use the constructions and notations of
the proof of Theorem \ref{thm.main} applied to $X$ and $\beta$.
In particular, $B$ denotes a topological complex rank 2 vector bundle
such that
$V\simeq \bP(B)$,
and hence
$$H^*(V)=H^*(X)[x]/(x^2-c_1(B)x+c_2(B)).$$

Consider the following claim:
the Picard number of $V$ is $r+1$ and the $\Q$-span of
the Neron-Severi group of $V$ is the $\Q$-span of $N$ together with
the element $2x+\gamma$, for some $\gamma\in H^2(X)$.
Given this, then, in
$$H^4(V)=H^4(X) \oplus xN
\oplus \langle x\beta\rangle \oplus xT,$$
we have
$$H^4(V)\cap H^{2,2}(V)=H^4(X) \oplus xN,$$
since for any $\delta\in N$,
$2(x\delta)=(2x+\gamma) \delta - \gamma \delta\in H^{2,2}(V)$.

Recall that we have
$f^*B\simeq E\otimes L$.
Now $x^2=x c_1(B) - c_2(B)=x c_1(E)+2xc_1(L)-c_2(B)$,
and by (\ref{c1class}),
$2c_1(L)\in\beta+2H^2(X)$, so the image of $x$ under the
obstruction map (\ref{obstr}) is nonzero.
The obstruction map for $V$ is thus completely determined by its
vanishing on elements of $H^2(X)$ and its nonvanishing on $x$.

It remains to verify the claim.
$\Pic({\mathcal G}\times_XV)$ is an extension of
$\Pic({\mathcal G})$ by the free group generated by $c_1(\cO_{\bP(E)}(1))$.
So $\Pic(V)$ is an extension of $\Pic(X)$ by some class which pulls back
to $2c_1(\cO_{\bP(E)}(1))$ on ${\mathcal G}\times_XV$.
But $2c_1(\cO_{\bP(E)}(1))$ differs from $2c_1(\cO_{\bP(f^*B)}(1))$
(the pullback of $2x$) by $2c_1(L)$, which lies in $H^2(X)$.
So, indeed, the Picard number of $V$ is $r+1$, and
for some $\gamma\in H^2(X)$,
$2x+\gamma$ is algebraic on $V$.
\end{proof}


\begin{thebibliography}{12}
\bibitem[1]{D-M} P. Deligne, D. Mumford, The irreducibility of
the space of curves of given genus, Publ. Math. IHES
{\bf 36} (1969), 75--109
\bibitem[2]{E-G} D. Edidin, W. Graham, Equivariant intersection
theory, Invent. Math. {\bf 131} (1998), 595--634
\bibitem[3]{FourAuthors} D. Edidin, B. Hassett, A. Kresch, A. Vistoli,
Brauer groups and quotient stacks, preprint (2000)
\bibitem[4]{E-L} R. Elman, T.Y. Lam,
Quadratic forms and the $u$-invariant, II,
Invent. Math. {\bf 21} (1973), 125-137
\bibitem[5]{Gr} A. Grothendieck, Le groupe de Brauer I, II, III,
in Dix Expos\'es sur la Cohomologie des Sch\'emas, pp. 46--188,
Adv. Stud. Pure Math. {\bf 3}, North-Holland, Amsterdam, 1968
\bibitem[6]{Kr} A. Kresch, Cycle groups for Artin stacks,
Invent. Math. {\bf 138} (1999), 495--536
\bibitem[7]{L-MB} G. Laumon, L. Moret-Bailly, Champs
Alg\'ebriques, Springer-Verlag, Berlin, 2000
\bibitem[8]{Me} A.S. Merkurev, On the norm residue symbol of
degree $2$, Dokl. Akad. Nauk SSSR {\bf 261} (1981), 542--547
\bibitem[9]{Milne} J. Milne, Etale Cohomology,
Princeton U. Press, Princeton, 1980
\bibitem[10]{Milnor} J. Milnor, The geometric realization of a
semi-simplicial complex, Ann. of Math. {\bf 65} (1957), 357--362
\bibitem[11]{Sal} D. Saltman, Noether's problem over an algebraically
closed field, Invent. Math. {\bf 77} (1984), 71--84
\bibitem[12]{Sar} V.G. Sarkisov, On conic bundle structures,
Izv. Akad. Nauk SSSR Ser. Math. {\bf 46} (1982), 371--408, 432
\end{thebibliography}
\end{document}